\documentclass[12pt,reqno]{amsart}
\usepackage{amsmath, amssymb, verbatim, comment, color, mathtools, accents}
\usepackage[margin=3cm]{geometry}

\makeatletter
\def\widebar{\accentset{{\cc@style\underline{\mskip10mu}}}}
\makeatother

\newcommand{\C}{\mathbb{C}}

\newcommand{\R}{\mathbb{R}}

\newcommand{\g}{\gamma}
\renewcommand{\d}{\delta}
\newcommand{\e}{\varepsilon}

\renewcommand{\t}{\theta}

\newcommand{\na}{\nabla}

\renewcommand{\o}{\omega}

\newcommand{\1}{({\rm I})}
\newcommand{\2}{({\rm II})}
\newcommand{\3}{({\rm III})}

\renewcommand{\i}{\sqrt{-1}}

\newcommand{\ddt}{\frac{d}{dt}}
\newcommand{\ddr}{\frac{d}{dr}}

\newcommand{\ie}{{\rm i.e.\ }} 

\DeclareMathOperator{\Ric}{Ric}

\DeclareMathOperator{\Tr}{Tr}

\newcommand{\lang}{\langle}
\newcommand{\rang}{\rangle}
\newcommand{\dmu}{d\mu}

\numberwithin{equation}{section}       

\theoremstyle{plain}
\newtheorem{theorem}{Theorem}[section]

\newtheorem{lemma}[theorem]{Lemma}

\theoremstyle{definition}

\newtheorem{example}[theorem]{Example}

\theoremstyle{remark}

\newtheorem*{acknowledgement}{\bf{Acknowledgment}} 

\title[Non-existence of eternal solutions]{Non-existence of eternal solutions to Lagrangian mean curvature flow with non-negative Ricci curvature} 
\date{\today}

\author[K. Kunikawa]{Keita Kunikawa}
\address{Advanced Institute for Materials Research \\
Tohoku University \\
Katahira \\
Aoba-ku \\
Sendai 980-8578 \\
Japan}
\email{keita.kunikawa.e2@tohoku.ac.jp}

\begin{document}

\begin{abstract}In this paper, we derive a mean curvature estimate for eternal solutions of uniformly almost calibrated Lagrangian mean curvature flow with non-negative Ricci curvature in the complex Euclidean space. As a consequence, we show a non-existence result for such eternal solutions. 
\end{abstract}

\maketitle
\section{Introduction}
Eternal solutions are solutions of mean curvature flow defined on a whole time interval $(-\infty, \infty)$. Since mean curvature flow is a quasi-linear parabolic equation, in general, we cannot solve the flow backward in time. Hence eternal solutions are very special and have interesting properties if they exist. Translating solutions are solutions of the mean curvature flow which have translation invariance under the mean curvature flow, and they are typical examples of eternal solutions. 

On the other hand, eternal solutions are known as type-II singularity models of mean curvature flow, namely, they arise after the parabolic rescaling at type-II singularities. In codimension one, Hamilton \cite{Ham95} showed that any complete convex eternal solution whose curvature attains its maximum at a point in the space-time must be a translating solution. However the situation is still not clear in higher codimensions. Therefore, in this paper, we consider not only translating solutions, but also eternal solutions to study type-II singularities. 

Besides above facts, Smoczyk \cite{Smo00} proved that there do not exist compact zero-Maslov class type-I singularities along the Lagrangian mean curvature flow. Later this result was extended to the complete case by Chen-Li \cite{CL01} and Wang \cite{Wan01} for the almost calibrated Lagrangian mean curvature flow, and then even later Neves \cite{Nev07} to the complete zero-Maslov class case. In these cases, therefore, we essentially need to investigate type-II singularities and eternal solutions (including translating solutions) as the singularity models. 

Our aim in this paper is to show a non-existence result for noncompact complete eternal solutions of almost calibrated Lagrangian mean curvature flow with an additional curvature condition. In this direction, we know some non-existence results by Han-Sun \cite{HS10}, Neves-Tian \cite{NT13} and Sun \cite{Sun13} for translating solutions. In \cite{HS10}, Han-Sun showed a non-existence result for almost calibrated Lagrangian translating solitons with non-negative sectional curvatures. We generalize their theorem \cite{HS10} to the class of almost calibrated Lagrangian eternal solutions with non-negative Ricci curvatures. 
\begin{theorem}\label{main} Let $\{L_t\}$ be a complete eternal solution to Lagrangian mean curvature flow in $\C^n$. Assume that $\{L_t\}$ is uniformly almost calibrated and 
\begin{align*}
\Ric_{L_t}\geq 0, \ \ \forall t\in(-\infty, \infty). 
\end{align*}
Then each $L_t$ must be totally geodesic. 
\end{theorem}

In general, the blow-up limit of a type-II singularity must be non-flat. Therefore, as a consequence of our result, any eternal solution with non-negative Ricci curvature cannot arise as a blow-up limit of uniformly almost calibrated Lagrangian mean curvature flow. 

Our result is a parabolic version of the theorem by Han-Sun \cite{HS10}. In order to show the parabolic version of the curvature estimate, we adopt techniques developed by Souplet-Zhang \cite{SZ06} and Wang \cite{Wan11}. Thanks to the parabolic curvature estimate, we do not need to assume the translation symmetry which was needed in \cite{HS10}. 

On the other hand, there are still some problems left in our paper. The assumption on Ricci curvature might be too strong to show the non-existence because of the following reasons: 
\begin{enumerate}
\item There do not exist any complete Lagrangian submanifolds with a strictly positive lower bound on the Ricci curvatures in Euclidean space. By Myers' theorem such a manifold would be compact with vanishing first Betti number which implies that the Maslov form (mean curvature form) is exact. Then by the Gauss equation, the Ricci curvatures must be non-positive in at least two points. In the author's knowledge, there is few source of examples for complete Lagrangians with non-negative Ricci curvatures. 
\item In general, non-negativity of Ricci curvature might not be preserved under the mean curvature flow. Therefore the assumption on the non-negativity of the Ricci curvature during the flow further narrows the available examples.  
\end{enumerate}
So far, we do not know whether the non-negativity assumption on Ricci curvature is essential or not.  

\begin{acknowledgement}
	The author is supported by Grant-in-Aid for JSPS Fellows Number 16J01498. During the preparation of this paper the author has stayed at the Max Planck Institute for Mathematics in the Sciences, Leipzig. The author is grateful to J\"{u}rgen Jost for his hospitality and his interest. Reiko Miyaoka also gave the author helpful comments in private seminars. 
\end{acknowledgement}

\section{Preliminaries}
\subsection{Lagrangian submanifold in $\mathbb{C}^n$}
Let $J$ and $\o(\cdot, \cdot)=\lang J\cdot, \cdot \rang$, respectively, be the standard complex structure on $\C^n$ and the standard symplectic form on $\C^n$. We consider the holomorphic $n$-form on $\C^n$ given by
\begin{align*}
\Omega:=dz^1 \wedge \cdots \wedge dz^n, 
\end{align*}
where $z^j=x^j+\i y^j\; (1\leq j \leq n)$ are the complex coordinates of $\C^n\cong\R^{2n}$. 

A smooth submanifold $L \subset (\C^n, \o)$ is said to be \emph{Lagrangian} if $\dim_\R=n$ and $\o|_L=0$. It is well known that the Lagrangian condition implies 
\begin{align*}
\Omega|_L=e^{\i\t}\dmu_L, 
\end{align*}
where $\dmu_L$ denotes the standard volume form of $L$ and $\t$ is some multivalued function called \emph{Lagrangian angle}. When the Lagrangian angle is single valued, $L$ is said to be in \textit{zero-Maslov} class, and moreover if $\cos \t >0$, 
then $L$ is said to be \textit{almost calibrated}. Especially, if there exists a positive constant $\d > 0$ such that 
\begin{align*}
\cos \t \geq \d > 0, 
\end{align*}
then we call $L$ \textit{strictly almost calibrated}. 

We denote by $A$ the second fundamental form, and $H:=\Tr A$ the mean curvature vector of the submanifold $L\subset \C^n$. Now, by the Gauss equation, we can write the scalar curvature $\mathrm{Scal}_L$ of $L$ as follows:  
\begin{align*}
\mathrm{Scal}_L=|H|^2-|A|^2. 
\end{align*}
Hence non-negative Ricci curvature on $L$ implies 
\begin{align}\label{AH}
|A|^2\leq|H|^2.
\end{align}
Note also that the important relation between the Lagrangian angle $\t$ and the mean curvature vector $H$ is given by 
\begin{align}
H=J\na\t. 
\end{align}
As an immediate consequence of this relation, we easily obtain the following: 
\begin{align}\label{cosH}
|\na(\cos\t)|^2\leq|H|^2. 
\end{align}

\subsection{Lagrangian mean curvature flow}
Next we consider a deformation of a Lagrangian submanifold $L^n \subset \C^n$. Let $F_0:L^n \to \mathbb{C}^{n}$ be a Lagrangian immersion. The \emph{mean curvature flow} is a one parameter family of smooth immersions $F:L^n\times[0, T_{\mathrm{max}})\rightarrow \mathbb{C}^{n}$ which satisfies the following:
\begin{align*}
\begin{cases}
\dfrac{d}{dt}F(p, t)=H(p, t),\quad p\in L^n, t\geq 0, \\
F(\cdot, 0)=F_0. 
\end{cases}
\end{align*} 
We denote by $L_t:=F_t(L):=F(L, t)$ a time slice of the flow $\{L_t\}_{t\geq 0}$. Note that mean curvature flow can be considered in more general situation. The submanifold is not necessarily Lagrangian and the ambient space is not necessarily Euclidean. However, it is shown by Smoczyk in \cite{Smo96} that the Lagrangian condition is preserved under mean curvature flow in K\"{a}hler--Einstein manifolds.  In this case we call the flow \textit{Lagrangian mean curvature flow}. It is also proved by Smoczyk \cite{Smo96, Smo99} that the zero-Maslov condition is preserved under the Lagrangian mean curvature flow, and the Lagrangian angle $\theta$ satisfies 
\begin{align*}
\ddt \t=\Delta\t.   
\end{align*}
These facts imply  
\begin{align}\label{evcos}
\bigg(\ddt-\Delta \bigg)\cos\t=|H|^2\cos\t.  
\end{align}
Therefore, by the parabolic maximum principle, the almost calibrated condition is also preserved under the Lagrangian mean curvature flow in Calabi--Yau manifolds. We say that a mean curvature flow $\{L_t\}$ is \textit{uniformly almost calibrated} if each time slice $L_t$ is strictly almost calibrated by a time independent constant $\d>0$. 

We note here another important evolution equation which was derived by Smoczyk in \cite{Smo00}: 
\begin{align}\label{evH}
\bigg(\ddt-\Delta\bigg)|H|^2=-2|\na^{\perp} H|^2+2|P|^2,  
\end{align}
where $\na^{\perp}$ denotes the connection on the normal bundle of $L$ and $P(\cdot, \cdot):=\lang A(\cdot, \cdot), H \rang$, therefore we have 
\begin{align*}
|P|^2\leq|A|^2|H|^2. 
\end{align*} 

\subsection{Eternal solutions}
A solution of the mean curvature flow is called \emph{eternal} if it is defined for all $t\in (-\infty, \infty)$. This kind of solution arises as a blow-up limit at a type II singularity of mean curvature flow. Note that any eternal solution in $\C^n$ must be noncompact, otherwise it has a finite time singularity. Thus, our attention is only in noncompact case. We call a solution \emph{Lagrangian eternal solution} if it is eternal and each time-slice is Lagrangian. All minimal submanifolds are stationary solutions of the mean curvature flow, thus they are trivial examples of eternal solutions. 

A submanifold $L^n \subset \C^n$ is called \emph{translating soliton} (or \emph{translator} for short) if there exists a nonzero constant vector $V \in \C^n$ so that 
\begin{align}\label{translating soliton}
H=V^\perp, 
\end{align}
where $V^\perp$ is the normal component of $V$ to the submanifold $L^n$. Let $F_0:L \to \C^n$ be a translating soliton and consider the following deformation of $L_0$: 
\begin{align}\label{translating solution}
F_t:=F_0+tV.  
\end{align}
It is clear that this deformation $L_t=F_t(L)$ is merely a parallel translation of the initial submanifold $L_0\subset\C^n$ in the constant direction of $V$, so that it is defined for all $t \in (-\infty, \infty)$. Furthermore we can show that \eqref{translating solution} satisfies the mean curvature flow. Conversely, it is known that if the solution of mean curvature flow can be written as the form \eqref{translating solution} for some submanifold $L_0$ and a constant vector $V\in\mathbb{C}^n$, then $L_0\subset\C^n$ must satisfy \eqref{translating soliton}, namely, a translating soliton. A solution of the mean curvature flow which has the form \eqref{translating solution} is called \textit{translating solution}. Of course, translating solutions are eternal solutions. By the translation symmetry, the study of a translating solution (a parabolic equation) is reduced to the study of a translating soliton (an elliptic equation) which is each time slice of the translating solution. 

\begin{example}
The curve in $\C \cong \R^2$ defined by 
\begin{align*}
\g_t:=\{(-\log\cos y+t, y)|-\pi/2<y<\pi/2\}, \ t\in (-\infty, \infty)
\end{align*}
is called the \emph{grim reaper} and it is a translating solution with the direction $V=(1, 0)$. It is known that the only translating solution in the plane is the grim reaper. 
\end{example}

\begin{example}
There is an eternal solution in $\mathbb{C}\cong\R^2$ given by 
\begin{align*}
\xi_t:=\{(x(t), y(t))\; |\; \sinh x(t)=e^{-t} \cos y(t)\}, \ t\in (-\infty, \infty), 
\end{align*} 
which is not a translating solution. It is called the \emph{hairclip} and it looks like infinitely many grim reapers, alternating between translating up and translating down for $t \to -\infty$, and converges to a vertical line as $t \to \infty$ (see for instance \cite{Has16}). 
\end{example} 

\begin{example}
Let $c_t\subset \C \cong \R^2$ be an eternal solution. Then a Lagrangian surface given by 
\begin{align*}
L_t:=c_t\times \mathbb{R} \subset \mathbb{C} \times \mathbb{C}=\C^2
\end{align*}
is also a (trivial) Lagrangian eternal solution in $\C^2$. 
\end{example}

Note that any curve in $\mathbb{C}$ is trivially Lagrangian and Ricci-flat.  In this case the Lagrangian angle $\t$ is the angle between the tangent vector of the curve and $x$-axis. Hence, above three eternal solutions give examples of Ricci-flat almost calibrated Lagrangian eternal solutions. However they all satisfy 
\begin{align*}
\inf_{L_t}\cos\theta=0, \quad \forall t\in(-\infty, \infty).  
\end{align*}
Therefore, the assumption $\cos\t \geq \d >0$ is needed in our main theorem. 

\section{Mean curvature estimate}
In this section we show a mean curvature estimate. As a consequence, we prove Theorem \ref{main}.  Our calculation is really similar to the gradient estimates in \cite{SZ06} and \cite{Wan11} for harmonic map heat flow on complete manifolds. Using their technique, the author showed a parabolic Bernstein type theorem for graphic eternal solutions of mean curvature flow (in codimension one) with bounded slope in \cite{Kun16}.  For graphic eternal solutions, the key point is to use the quantities $|A|^2$ and the slope of the graph to show the parabolic curvature estimate. In our current case, that is, for almost calibrated Lagrangian eternal solution, it turns out that $|H|^2$ and Lagrangian angle $\cos\t$ match in order for the parabolic curvature estimate. 

\subsection{Mean curvature estimate on a cylindrical domain}
Let $\{L_t\}_{t \in [-T, T]}$ be a complete solution to almost calibrated Lagrangian mean curvature flow with $L_0 = L$, a complete Lagrangian submanifold in $\C^n$.  We denote by $r(p):=d(o, p)$ the intrinsic distance on $L$ from a fixed point $o \in L$ to a point $p \in L$. Of course if $p\in L$ is not in cut locus of $o \in L$, the distance function $r$ is differentiable and $|\na r| = 1$. Let $D_R=D_R(o)$ be a closed intrinsic distance ball with center $o$ and radius $R > 0$. Moreover, take a space-time cylindrical domain $D_{R, T}=D_{R, T}(o):=D_R(o)\times[-T, T]$.  Define the function $\varphi(p, t):=1-\cos\t(p, t)$ on $\{L_t\}_t$. Now our mean curvature estimate on the cylindrical domain is given as follows. 
\begin{theorem}\label{curvestimate}
Let $F:L^{n}\times[-T, T] \to \C^n$ be a complete solution to almost calibrated Lagrangian mean curvature flow which has non-negative Ricci curvature for all time.  Assume further that there exists a positive constant $\d > 0$ such that $\cos\t(p, t)\geq \d > 0$ for any point in $L^n\times[-T, T]$. Then there exists a constant $C$ which is independent of $R$ and $T$ such that
\begin{align*}
\sup_{D_{R/2, T/2}}\frac{|H|}{b-\varphi}\leq C\Big(\frac{1}{R}+\frac{1}{\sqrt{R}}+\frac{1}{\sqrt{T}}\Big), 
\end{align*}
where $b$ is a constant satisfying $\sup_{L_t} \varphi\leq 1-\d < b < 1$ for $t\in[-T, T]$. 
\end{theorem}

For a complete eternal solution $\{L_t\}$ satisfying the assumption in Theorem \ref{main}, there exists a global constant $b$ such that $\sup_{L_t} \varphi\leq 1-\d < b < 1$, for all $t\in(-\infty, \infty)$. Take $R \to  \infty$ and $T \to \infty$ in Theorem \ref{curvestimate}, then it follows $|H|\to 0$. Finally, we  obtain $|A|\equiv 0$ under the assumption $|A|^2\leq|H|^2$. Therefore we obtain Theorem \ref{main}. 

\subsection{Proof of the mean curvature estimate}
We prove Theorem \ref{curvestimate} here. 
\begin{proof}
We assume that $\{L_t\}$ is uniformly almost calibrated, namely, there exists a positive constant $\d > 0$ so that 
\begin{align*}
\cos\t \geq\d >0, \ \ \forall t\in[-T, T]. 
\end{align*}
Set $\varphi:=1-\cos\t$. Then since $\sup \varphi\leq1-\d$, there exists a constant $b$ such that
\begin{align*}
\sup_{L_t} \varphi\leq 1-\d < b < 1, \ \ \ \forall t\in[-T, T]. 
\end{align*}
We also assume that $L_t$ has non-negative Ricci curvature for each $t\in[-T, T]$.  Therefore the inequality \eqref{AH} always holds, that is, $|A|^2\leq|H|^2$ for all $t\in[-T, T]$. Combining \eqref{AH} and \eqref{evH} we obtain
\begin{align}\label{evH2}
\bigg(\ddt-\Delta\bigg)|H|^2\leq-2|\na^{\perp} H|^2+2|H|^4. 
\end{align}
Define the function on $\{L_t\}_{t\in[-T, T]}$ by
\begin{align*}
\phi=\frac{|H|^2}{(b-\varphi)^2}. 
\end{align*}
A direct calculation shows that 
\begin{align}
\na\phi=\frac{\na|H|^2}{(b-\varphi)^2}+\frac{2|H|^2\na\varphi}{(b-\varphi)^3}. 
\end{align}
Similarly we can compute 
\begin{align}\label{me}
\Delta\phi=\frac{\Delta|H|^2}{(b-\varphi)^2}+\frac{4\lang\na\varphi, \na|H|^2\rang}{(b-\varphi)^3}
+\frac{2|H|^2\Delta\varphi}{(b-\varphi)^3}+\frac{6|\na\varphi|^2|H|^2}{(b-\varphi)^4}. 
\end{align}
On the other hand, the time derivative of $\phi$ is given by
\begin{align}\label{meee}
\ddt\phi=\frac{\ddt|H|^2}{(b-\varphi)^2}-\frac{2|H|^2\ddt\cos\t}{(b-\varphi)^3}. 
\end{align}
Subtracting \eqref{meee} from \eqref{me},  we obtain
\begin{align*}
\bigg(\Delta-\ddt\bigg)\phi=\frac{\Delta|H|^2-\ddt|H|^2}{(b-\varphi)^2}+&\frac{4\lang \na\varphi, \na|H|^2\rang}{(b-\varphi)^3}\\+&\frac{2|H|^2\ddt\cos\t-2|H|^2\Delta\cos\t}{(b-\varphi)^3}+\frac{6|\na\varphi|^2|H|^2}{(b-\varphi)^4}. 
\end{align*}
By using \eqref{evcos} and \eqref{evH2}, we compute 
\begin{align*}
\bigg(\Delta-\ddt \bigg)\phi\geq\frac{2|\na^{\perp}H|^2-2|H|^4}{(b-\varphi)^2}+&\frac{4\lang \na\varphi, \na|H|^2\rang}{(b-\varphi)^3}\\ +&\frac{2|H|^4\cos\t}{(b-\varphi)^3}+\frac{6|\na\varphi|^2|H|^2}{(b-\varphi)^4}. 
\end{align*}
Note that the following relations hold, 
\begin{align}
\frac{2|\na^{\perp} H|^2}{(b-\varphi)^2}+\frac{2|\na\varphi|^2|H|^2}{(b-\varphi)^4}&\geq\frac{4|\na^{\perp} H||\na\varphi||H|}{(b-\varphi)^3}, \label{mei}\\
\frac{2\lang\na|H|^2, \na\varphi\rang}{(b-\varphi)^3}+\frac{4|H|^2|\na\varphi|^2}{(b-\varphi)^4}&=
\frac{2\lang\na\varphi, \na\phi\rang}{(b-\varphi)}. \label{meimei}
\end{align}
By using \eqref{mei} and \eqref{meimei} with the Cauchy-Schwarz inequality and Kato's inequality $|\nabla|H||^2\leq|\nabla^{\perp}H|^2$, we finally obtain
\begin{align}\label{meimeimei}
\bigg(\Delta-\ddt \bigg)\phi\geq\frac{2(1-b)|H|^4}{(b-\varphi)^3}+\frac{2\lang\na\varphi, \na\phi\rang}{b-\varphi}. 
\end{align}

Now we take a smooth function $\eta(r, t):\mathbb{R}\times \mathbb{R}\rightarrow\mathbb{R}$ supported on $[-R, R]\times[-T, T]$ which has the following properties: 
\begin{align*}
(\mathrm{i}) &\ \eta(r, t)\equiv 1 \textrm{ on } [-R/2, R/2]\times[-T/2, T/2] \textrm{ and } 0\leq\eta\leq1, \\
(\mathrm{ii}) &\ \eta(r, t) \textrm{ is decreasing if } r\geq 0 \quad \ie \ddr\eta\leq0, \\
(\mathrm{iii}) &\ \Big|\ddr \eta\Big|\eta^{-a}\leq C_a R^{-1}\quad \text{and} \quad \Big|\frac{d^2}{dr^2}\eta\Big|\eta^{-a}\leq C_a R^{-2} \quad \textrm{for } 0<a<1, \\  
(\mathrm{iv}) &\  \Big|\ddt\eta\Big|\eta^{-a}\leq C_a T^{-1} \quad \textrm{for } 0<a<1. 
\end{align*}
Such a function can be given by a standard way (see \cite{LY86},  \cite{SZ06} or \cite{Wan11}). We use a cut-off function supported on $D_{R, T}$ given by
$\psi(p, t):=\eta(r(p), t)$. 

Let $\Phi:=-2\na \varphi/(b-\varphi)$. By using \eqref{meimeimei} we can calculate 
\begin{align}\label{meei}
\Delta(\psi &\phi)+\lang \Phi, \na(\psi \phi)\rang-2\Big\lang \frac{\na\psi}{\psi}, \na(\psi \phi)\Big\rang-\ddt(\psi \phi) \\
&=\psi \bigg(\Delta\phi-\ddt\phi\bigg)+\phi\bigg(\Delta\psi-\ddt\psi\bigg)+\lang\psi \Phi, \na\phi\rang+\lang\phi \Phi, \na\psi\rang-2\frac{|\na\psi|^2}{\psi}\phi \nonumber\\
&\geq 2\psi\frac{(1-b)|H|^4}{(b-\varphi)^3}+\phi\bigg(\Delta\psi-\ddt\psi\bigg)-2\frac{\lang\na\varphi, \na\psi\rang}{b-\varphi}\phi-2\frac{|\na\psi|^2}{\psi}\phi. \nonumber
\end{align}
Since $D_{R, T}$ is compact, $\psi\phi$ attains its maximum at some point $(p_0, t_0)$ in $D_{R, T}$. We may assume that such a point $p_0\in L$ is not in the cut locus of $o\in L$ (see \cite{CY75}, \cite{SZ06} or \cite{Wan11}). At the point $(p_0, t_0)\in D_{R, T}$, we have 
\begin{align*}
\na(\psi \phi)=0,\quad  \Delta(\psi \phi)\leq0,\quad   \ddt(\psi \phi)\geq0. 
\end{align*}
Hence by \eqref{meei}, we obtain
\begin{align}
2\psi(1-b)\frac{|H|^4}{(b-\varphi)^3}&\leq 2\phi\frac{\lang\na\varphi, \na\psi\rang}{b-\varphi}+2\phi\frac{|\na\psi|^2}{\psi}+\phi(\ddt\psi-\Delta\psi) \nonumber\\
&=:\1+\2+\3. \nonumber
\end{align}
Note that the following holds: 
\begin{equation*}
|\nabla\psi|^2=\bigg|\ddr \eta\bigg|^2|\nabla r|^2\leq \bigg|\ddr\eta\bigg|^2. 
\end{equation*}
Recall that the following Young's inequality.  
\begin{lemma}[Young's inequality]
For any $a, b > 0$ and any $\e > 0$, we have
\begin{align*}
ab\leq\e a^p+\e^{-\frac{q}{p}}b^q, 
\end{align*}
where $1<p, q<\infty$ and $1/p+1/q=1$. 
\end{lemma}
By using \eqref{cosH}, Young's inequality and the property of $\eta$, we can estimate $\1$ as follows: 
\begin{align}\label{I}
\1&\leq2\phi\frac{|\na\varphi|}{b-\varphi}|\na\psi|\leq2\phi\frac{|H|}{b-\varphi}|\na\psi|=2\phi^{\frac{3}{2}}|\na\psi|\\
&\leq\frac{\e}{4}\psi\phi^2+\frac{C(\e)|\na\psi|^4}{\psi^3}
\leq\frac{\e}{4}\psi\phi^2+\frac{C(\e)|\ddr\eta|^4}{\psi^3}\nonumber\\
&\leq\frac{\e}{4}\psi\phi^2+\frac{C(\e)}{R^4}, \nonumber
\end{align}
where $\e>0$ is an arbitrary constant which is determined later and $C(\e)$ are constants depending only on $\e$. Similarly, we can compute by using Young's inequality and the property of $\eta$, 
\begin{align}\label{II}
\2=2\phi\frac{|\na\psi|^2}{\psi}\leq\frac{\e}{4}\psi\phi^2+C(\e)\frac{|\na\psi|^4}{\psi^3}
\leq\frac{\e}{4}\psi\phi^2+\frac{C(\e)}{R^4}. 
\end{align}
Since the Ricci curvature of $L$ is non-negative, by the Laplacian comparison theorem we have
\begin{align}\label{lapcomp}
\Delta r\leq\frac{n-1}{r}. 
\end{align}
Using \eqref{lapcomp} and $\ddr \eta\leq 0$, we have
\begin{align*}
\Delta\psi=\Delta r\bigg(\ddr \eta\bigg)+|\na r|^2\bigg(\frac{d^2}{dr^2}\eta\bigg)\geq\frac{n-1}{r}\bigg(\ddr \eta\bigg)-\bigg|\frac{d^2}{dr^2}\eta\bigg|. 
\end{align*}
Hence we obtain for the second term of $\3$ in the same way as above, 
\begin{align}\label{III1}
-\phi\Delta\psi&\leq \phi\bigg|\frac{d^2}{dr^2}\eta\bigg|+\frac{(n-1)|\ddr \eta|\phi}{r}\\
&\leq \phi\bigg|\frac{d^2}{dr^2}\eta\bigg|+\frac{2(n-1)|\ddr\eta|\phi}{R} \nonumber\\
&\leq \frac{\e}{4}\psi\phi^2+C(\e, n)\bigg(\frac{1}{R^4}+\frac{1}{R^2}\bigg),  \nonumber
\end{align}
where $C(\e, n)$ is a constant depending only on $\e$ and $n$. 
(Note that we may assume $R/2\leq r$ for the second inequality since $\ddr\eta\equiv0$ for $r\leq R/2$.)

As for the first term of $\3$, we have
\begin{align}\label{III2}
\phi\bigg(\ddt\psi\bigg)&\leq\phi\bigg|\ddt \eta\bigg|
\leq\frac{\e}{4}\psi\phi^2+C(\e)\frac{|\ddt\eta|^2}{\psi}\\
&\leq \frac{\e}{4}\psi\phi^2+\frac{C(\e)}{T^2}. \nonumber
\end{align}
From \eqref{I}, \eqref{II}, \eqref{III1} and \eqref{III2}, we finally obtain
\begin{align*}
2(1-b)(b-\varphi)\psi\phi^2\leq\e\psi\phi^2+C(\e, n)\bigg(\frac{1}{R^4}+\frac{1}{R^2}+\frac{1}{T^2}\bigg). 
\end{align*}
Since $\e > 0$ is arbitrary we can take a sufficiently small $\e$ such that 
\begin{align*}
2(1-b)(b-\varphi)-\e>0. 
\end{align*}
Then we have
\begin{align*}
(\psi\phi)^2\leq\psi\phi^2\leq C\bigg(\frac{1}{R^4}+\frac{1}{R^2}+\frac{1}{T^2}\bigg). 
\end{align*}
Since $\psi\equiv 1$ on $D_{R/2, T/2}$, 
\begin{align*}
\sup_{D_{R/2, T/2}}\frac{|H|}{b-\varphi}\leq C\bigg(\frac{1}{R}+\frac{1}{\sqrt{R}}+\frac{1}{\sqrt{T}}\bigg). 
\end{align*}
This completes the proof of Theorem \ref{curvestimate}. 
\end{proof}


\end{document}